\documentstyle[12pt]{article}

\newtheorem {theorem}{Theorem}
\newtheorem {lemma}{Lemma}
\newtheorem {corollary}{Corollary}

\newcounter {rem}

\newcounter {exa}
\newenvironment {example}
{\refstepcounter{exa}\par\medskip\noindent{\bf Example
\arabic{exa}\,\,} }{\par}
\newcounter {def}

\newenvironment {Remark} {\par\medskip\noindent{\bf
Remark\,\,}} {\par}

\newenvironment {proof} {\par\medskip\noindent{\bf
Proof\,\,}} {\qed\par}

\newenvironment {Definition} {\par\medskip\noindent{\bf
Definition\,\,}}{\par}

\def \Real {I\!\!R}

\def \Complex {\,\prime\mskip-2.5\thinmuskip C}

\def \lnorm#1\rnorm {\vphantom{#1}\left\|\smash{#1}\right\|}
\def \lmod#1\rmod {\vphantom{#1}\left|\smash{#1}\right|}

\newcommand \bydef {\stackrel{\mbox{\scriptsize def}}{=}}
\newcommand \qed {\,\rule[-.23ex]{1.6ex}{1.6ex}}
\newcommand \Planck {{h^{\mskip-4.5\thinmuskip{-}}}}

\renewcommand \phi {\varphi}
\renewcommand \rho {\varrho}

\author{A.A. Balinsky\thanks{E-mail :
balin@leeor.technion.ac.il}, Yu.M.  Burman\thanks{E-mail :
mat9304@technion.technion.ac.il}\\
{\normalsize Department of Mathematics}\\
{\normalsize Technion --- Israel Institute of Technology,}
\\
{\normalsize 32000, Haifa, Israel}}
\title{Quadratic Poisson brackets and Drinfel'd theory for
associative algebras}
\date{}
\begin{document}

\font\eux=eufm10
\font\eus=eufm8
\font\euss=eufm6
\newfam\eufam
\textfont\eufam=\eux
\scriptfont\eufam=\eus
\scriptscriptfont\eufam=\euss
\def\eu{\fam\eufam}

\newcommand \CC {C^\infty}
\newcommand \Mat {\mathop{\rm Mat}}
\newcommand \Symm {\mathop{\rm Symm}}
\newcommand {\SchB}[1]{[[#1,#1]]}
\newcommand {\ii} {{\bf i}}
\newcommand {\jj} {{\bf j}}
\newcommand {\kk} {{\bf k}}

\maketitle

\centerline{\bf Abstract}\bigskip

{\small Quadratic Poisson brackets on associative algebras are studied.
Such a bracket compatiblle with the multiplication is related
to a differentiation in tensor square of the underlying algebra.
Jacobi identity means that this differentiation satisfies
a classical Yang--Baxter equation. Corresponding Lie groups are
canonically equipped with a Poisson Lie structure. A way to
quantize such structures is suggested.}

\section{Introduction}
The main purpose of the present article is to construct an
analog of Drinfel'd theory (see \cite{DR1}) for associative
algebras. Poisson Lie brackets involved in the original version
of the theory are replaced by {\em quadratic} Poisson brackets
compatible with the algebra structure. Restriction of such
brackets to the group of invertible elements of the algebra
leads to a Poisson Lie structure (see Theorem \ref{Group}
below).

As usual, Poisson bracket $\{\cdot,\cdot\}$ on a smooth manifold
$M$ is understood as a Lie algebra structure on the space of
smooth functions $\CC(M)$ satisfying the Leibnitz identity
$\{fg,h\} = f\{g,h\} + \{f,h\}g$. If Jacobi identity is not
required one speaks about pre-Poisson brackets. A mapping $F:M_1
\to M_2$ of two manifolds equipped with Poisson brackets
$\{\cdot,\cdot\}_1$ and $\{\cdot,\cdot\}_2$, respectively, is
called {\em Poisson} if a canonical mapping $F^*: \CC(M_2) \to
\CC(M_1)$ is a Lie algebra homomorphism, i.e.
\begin{displaymath}
\{f\circ F,g \circ F\}_1 = \{f,g\}_2 \circ F.
\end{displaymath}

Suppose now that $M$ is equipped with a multiplication, i.e. a
mapping $* \colon M \times M \to M$. This immediately gives rise
to a co-multiplication (diagonal) $\Delta: \CC(M) \to \CC(M)
\otimes \CC(M)$ acting as
\begin{displaymath}
\Delta(f)(x,y) = f(x * y),
\end{displaymath}
here a usual identification of $\CC(M \times M)$ with $\CC(M)
\otimes \CC(M)$ is assumed.

A Poisson bracket is said to be {\em compatible} with the
multiplication if $* \colon M \times M \to M$ is a Poisson
mapping where $M \times M$ is equipped with a product Poisson
bracket. In other words, the following identity should be
satisfied:
\begin{equation}\label{DiagComp}
\Delta(\{f,g\}) = \{\Delta(f),\Delta(g)\}
\end{equation}
where
\begin{displaymath}
\{p \otimes q, r \otimes s\} \bydef \{p,r\} \otimes qs + pr
\otimes \{q,s\}.
\end{displaymath}

In particular, a bracket on a vector space compatible with
addition structure is exactly a Berezin--Lie one, i.e. a linear
bracket (speaking about linear, quadratic, etc., brackets we
will always mean that a bracket of two {\em linear} functions is
linear, quadratic, etc., respectively).

The next-simplest case is quadratic Poisson brackets (see
\cite{Skl}, \cite{BAL2}, \cite{KUP2}). Such a bracket may be
compatible only with a bilinear operation on the vector space,
i.e. with an algebra structure. One of the most important cases,
that of a full matrix algebra $\Mat(n,K)$ where $K = \Real$ or
$\Complex$, was investigated in detail in a series of works by
B.~ Kupershmidt \cite{KUP3,KUP4,KUP5}.

A natural question arises, whether it is possible to obtain a
quantization of these brackets. This can be done using quadratic
algebras (see \cite{Manin}) and will be a subject of the
forthcoming paper. Quantization of the corresponding Poisson Lie
groups may also be simpler than in general case (see e.g.
\cite{DR}) because such Poisson Lie groups bear a preferred
coordinate system (linear functions on the original algebra).

Authors are grateful to J. Donin, B. Kupershmidt, A. Stolin, and
N. Zobin for many stimulating discussions.

\section{Quadratic brackets compatible with algebra structure}
Let $A$ be an associative algebra, i.e. vector space with a
bilinear associative multiplication $*$. The symbol $A_L$ will
denote an {\em adjacent} Lie algebra, i.e. Lie algebra on the
same space with a commutator $[a,b] = ab - ba$. The dual space
$A^*$ is naturally embedded into functions algebra $\CC(A)$, as
well as all its symmetric powers $\Symm\bigl((A^*)^{\otimes
n}\bigr)$ are.  Again, Poisson bracket is understood as a Lie
algebra structure on $\CC(A)$, i.e. mapping $B \colon \CC(A)
\wedge \CC(A) \to \CC(A)$ satisfying Jacobi identity.

\begin{Definition}
A Poisson bracket $B$ is called {\em quadratic} if
\begin{displaymath}
B(A^* \wedge A^*) \subset \Symm\bigl((A^*)^{\otimes 2}\bigr).
\end{displaymath}
\end{Definition}
If $\{e_i\}$ is a basis in $A$ and $x^i$ are coordinate
functions, then quadratic Poisson bracket is given by
\begin{equation}\label{DefPoi}
\{x^i,x^j\} = c_{kl}^{ij}x^kx^l,
\end{equation}
a summation over repeated indices will be always assumed.

With $A$ being an algebra, its tensor square $A \otimes A$ can
also be given an algebra structure by a componentwise
multiplication. Then the set of symmetric tensors $\Symm(A
\otimes A) \subset A \otimes A$ is its subalgebra, while the set
of skew-symmetric tensors $A \wedge A \subset A \otimes A$ is a
$\Symm(A \otimes A)$-bimodule. A linear mapping $D: {\eu A} \to
V$ from algebra $\eu A$ to a $\eu A$-bimodule $V$ is called a
{\em differentiation} if it obeys the condition
\begin{equation}\label{DefDif}
D(p \cdot q) = p D(q) + D(p) q
\end{equation}
for all $p,\,q \in \eu A$.

\begin{theorem}\label{Diff}
Quadratic Poisson bracket $\delta^* \colon A^* \wedge A^* \to
\Symm\bigl((A^*)^{\otimes 2}\bigr)$ is compatible with the
algebra structure in $A$ if and only if its dual mapping
\begin{displaymath}
\delta \colon \Symm\bigl(A^{\otimes 2}\bigr) \to A \wedge A
\end{displaymath}
is a differentiation.
\end{theorem}

See \cite{BAL2} for proof. It does not make use of Jacobi
identity, and therefore applies for general pre-Poisson brackets
as well.

\begin{example}\cite{FRT}
Consider a Poisson bracket
\begin{displaymath}
\{x^i,x^j\} = x^ix^j \quad \mbox{for $i < j$}.
\end{displaymath}
It is compatible with the multiplication of the ``first column
algebra'' (algebra of $n \times n$-matrices having nonzero
elements only in the first column).
\end{example}

For a linear operator $P \colon V \otimes V \to V \otimes V$
symbol $P^{12}$ will mean operator $V \otimes V \otimes V \to V
\otimes V \otimes V$ acting as $P$ on the first and second
tensor component and as identity operator on the third one;
notations like $P^{13}$ have a similar meaning. A {\em Schouten
bracket} $\SchB{P}$ is then defined as
\begin{displaymath}
[P^{12},P^{13}] + [P^{12},P^{23}] + [P^{13},P^{23}].
\end{displaymath}
with $[\cdot,\cdot]$ meaning an ordinary commutator.

Let $\delta^*$ be a quadratic pre-Poisson bracket, and
$\widetilde \delta:  A \otimes A \to A \otimes A$ be an
arbitrary linear extension of the dual operator $\delta$.

\begin{theorem}\label{Schou}
Bracket $\delta^*$ is Poisson (i.e. satisfies Jacobi identity)
if and only if
\begin{equation}\label{ESchou}
\SchB{\widetilde \delta}(X) = 0
\end{equation}
for any fully symmetric tensor $X \in A \otimes A \otimes A$.
\end{theorem}
\begin{Remark}
In particular, it is asserted in Theorem that $\SchB{\widetilde
\delta}(X)$ for fully symmetric $X$ depends on $\delta$ only,
and not on its specific extension $\widetilde\delta$.
\end{Remark}

The proof is made by a straightforward (and rather tedious)
computation.

Note also that Theorem \ref{Schou} does not require bracket
$\delta^*$ to be compatible with algebra.

\section{Quadratic brackets and Poisson Lie groups}
Lie bialgebra structures (cocommutators) on a Lie algebra $L$
are in one-to-one correspondence with the {\em linear} Poisson
brackets compatible with a commutator in $L$ (see \cite{DR1}).
Quadratic Poisson brackets on associative unital algebras are
also directly connected with Lie bialgebras:

\begin{lemma}\label{BiAlg}
Let $A$ be an associative algebra with a unit $u$, and $\delta^*
\colon A^* \wedge A^* \to \Symm\bigl((A^*)^{\otimes 2}\bigr)$ be
a quadratic Poisson bracket compatible with it. Then the mapping
(comultiplication)
\begin{equation}\label{CoMult}
\Delta(x) = \delta(x \otimes u + u \otimes x)
\end{equation}
equips $A_L$ with a Lie bialgebra structure.
\end{lemma}

Proof: computation.

So, a differentiation $\delta \colon \Symm\bigl(A^{\otimes
2}\bigr) \to A \wedge A$ satisfying condition (\ref{ESchou})
defines two brackets compatible with $A_L$, quadratic $\delta^*$
and linear $\Delta^*$ with $\Delta$ given by (\ref{CoMult}).

\begin{theorem}\label{Compat}
These two brackets are compatible with one another, i.e. their
arbitrary linear combination $\alpha \delta^* + \beta \Delta^*$
is also a Poisson bracket compatible with $A_L$.
\end{theorem}
\begin{proof}
Let $u = u^i e_i$ be a unit of $A$ (as usual, $\{e_i\} \subset
A$ is an additive basis). Consider a change of coordinates
$f_t(x) = x - tu$ ($t \in \Real$ is a parameter). In new
coordinates the bracket $\delta^*$ becomes $\delta^* +
t\Delta^*$. This obviously means a compatibility.
\end{proof}

The following Theorem directly relates quadratic brackets to
Poisson Lie structures. Drinfel'd in \cite{DR1} shows that there
is a one-to-one correspondence between Poisson Lie groups and
Lie bialgebra structures.

\begin{theorem}\label{Group}
Let $A$ be an associative algebra with a unit $u$, and $\delta^*
\colon A^* \wedge A^* \to \Symm\bigl((A^*)^{\otimes 2}\bigr)$ be
a quadratic Poisson bracket compatible with it. Then the group
of invertible elements of $A$ is a Poisson Lie group with
respect to a restriction of the bracket $\delta^*$ to it. The
corresponding bialgebra is exactly as in Lemma \ref{BiAlg}.
\end{theorem}
\begin{proof}
The group is an open submanifold of $A$, so, it is apparently a
Poisson submanifold. Since the group inherits its multiplication
from the algebra, it is a Poisson Lie group. To obtain bialgebra
structure, on has to take a bracket of two coordinate function
in the vicinity of the unit $u$, and to take its differential at
the point $u$. Again, take a change of variables $f_t(x) = x -
tu$, and the bracket of cordinate functions will be $\delta^* +
t\Delta^*$; the differential at the point $u$ is exactly
$\Delta$.
\end{proof}

\begin{corollary}\label{QuadCoor}
Let $G$ be a Poisson Lie group whose Lie algebra $\eu G$ is
adjacent to an unital associative algebra structure. Then $G$
has a coordinate system in which its Lie Poisson bracket is
quadratic.
\end{corollary}

\section{Coboundary brackets}
Restrict now our consideration to a {\em coboundary} case when
differentiation $\delta$ is internal i.e. given by a formula
\begin{equation}\label{Commut}
\delta(x) = [r,x]
\end{equation}
with some $r \in A \wedge A$ (mapping (\ref{Commut}) is
obviously a differentiation).

If $r = \sum_i \alpha_i \otimes \beta_i \in A \otimes A$ then
$r^{12}$ is, by definition, $\sum_i \alpha_i \otimes \beta_i
\otimes 1 \in A \otimes A \otimes A$. Here $1$ is a formal unit
(a unit of an universal enveloping algebra of $A_L$). Then on
Schouten bracket $\SchB{r} \in A \otimes A \otimes A$ is defined
as
\begin{equation}\label{SchouEl}
[r^{12},r^{13}] + [r^{12},r^{23}] + [r^{13},r^{23}].
\end{equation}

\begin{theorem}\label{MYBE}
Differentiation $\delta$ satisfies condition (\ref{ESchou}) (and
therefore defines a Poisson bracket) if and only if $\SchB{r}$
commutes with any fully symmetric tensor $X \in A \otimes A
\otimes A$.
\end{theorem}

\begin{Remark}
If $A$ has a unit $1$ then hypothesis of Theorem \ref{MYBE} is
equivalent to the following property: $\SchB{r}$ commutes with
all the tensors of the type $1 \otimes 1 \otimes x + 1 \otimes x
\otimes 1 + x \otimes 1 \otimes 1$, for $x \in A$. Or, in other
words, $\SchB{r}$ should be $\rm Ad$-invariant with respect to
the Lie algebra $A_L$.
\end{Remark}

Theorem \ref{MYBE} is just Theorem \ref{Schou} for a coboundary
bracket.  Note also that for coboundary quadratic Poisson
brackets one can introduce concepts of triangular,
quasi-triangular, factorizable, etc., brackets, like those
introduced in the theory of coboundary Lie bialgebras (see
\cite{DR}).

\begin{corollary}\label{SubAlg}
Let $A$ be associative algebra, and $\eu G$ is a Lie subalgebra
of $A_L$. Then any element $r \in {\eu G} \wedge {\eu G}$
satisfying classical Yang--Baxter equation (see \cite{DR})
defines a quadratic Poisson bracket compatible with $A$.
\end{corollary}

\section{Examples}
Besides those described here, several important examples can be
found in the classical article \cite{STS1}.

\begin{example}\label{ab-ba}
Let $A$ be an associative algebra with a solvable $A_L$. By Ado
theorem any representation of $A_L$ has an eigenvector. Applying
this to a regular representation, one obtains that $A_L$
contains two elements, $a$ and $b$, with $[a,b] = sb$ for some
constant $s$. The element $r = a \otimes b - b \otimes a$
satisfies classical Yang--Baxter equation and thus, according to
Corollary \ref{SubAlg}, defines a quadratic Poisson bracket on
$A$.
\end{example}

\begin{example}\label{Quaternions}
Consider a body $\bf H$ of quaternions as a four-dimensional
$\Real$-algebra. Then any element $r \in {\bf H} \wedge {\bf H}$
of the type
\begin{equation}\label{rquatern}
r = a\,\ii \wedge \jj + b\,\ii \wedge \kk + c\,\jj \wedge \kk
\end{equation}
satisfies conditions of Theorem \ref{MYBE}. Corresponding
Poisson bracket on $\bf H$ is:
\begin{eqnarray}
\{x^1,x^2\} &=& x^2(bx^3 - ax^4) + c((x^3)^2 + (x^4)^2)
\label{RowFirst}\\
\{x^1,x^3\} &=& -x^3(cx^2 + ax^4) - b((x^2)^2 + (x^4)^2)\\
\{x^1,x^4\} &=& x^4(-cx^2 + bx^3) + a((x^2)^2 + (x^3)^2)\\
\{x^2,x^3\} &=& x^1(-bx^2 + cx^3)\\
\{x^2,x^4\} &=& -x^1(ax^2 + cx^4) \\
\{x^3,x^4\} &=& x^1(ax^3 - bx^4)\label{RowLast}
\end{eqnarray}
where $x^1,\,x^2,\,x^3,\,x^4$ is a basis in ${\bf H}^*$ dual to
$1,\,\ii,\,\jj,\,\kk \in {\bf H}$.

A Lie group corresponding to $\bf H$ is a multiplicative group
of all nonzero quaternions. It contains a subgroup $\{z \in {\bf
H}, \lnorm z\rnorm = 1\}$ isomorphic to $\mathop{\rm SU}(2)$. It
can be checked (see \cite{BAL2}) that this subgroup is Poisson
with respect to all the brackets
(\ref{RowFirst})--(\ref{RowLast}). Thus, $\mathop{\rm SU}(2)$
bears a three-parameter family of Poisson Lie structures.
\end{example}

\begin{example}\label{Nilpot}
Consider a Lie algebra $\eu G$ satisfying an identity $[[{\eu
G},{\eu G}], {\eu G}] = 0$ (i.e. $[[a,b],c] = 0$ for all
$a,\,b,\,c \in \eu G$) (a particular case here is Heisenberg
algebra, a three-dimensional Lie algebra with generators
$p,\,q$, and $z$, and relations $[p,q] = \Planck z$, and $[z,
\cdot] = 0$). Then $a * b = \frac{1}{2} [a,b]$ is an associative
operation with $[a,b] = a*b - b*a$.  An arbitrary element $r \in
{\eu G} \wedge {\eu G}$ commutes with all the symmetric tensors
from ${\eu G} \otimes {\eu G}$ and therefore defines a zero
Poisson bracket on $\eu G$.

Associative algebra $\eu G$ has no unit. Having it added,
consider a space $A = \langle 1\rangle \oplus {\eu G}$. Then the
Lie group $G$ corresponding to $\eu G$ is an affine subspace $1
+ {\eu G} \subset A$ (it is the usual implementation of $G$ via
exponents in the universal enveloping algebra of $\eu G$; the
thing is that $a^2 = 0$ for any $a \in {\eu G}$).  Fix a basis
$\{e_i\} \in {\eu G}$; it defines a natural coordinate system in
the group $G$. The above construction applied to $r = r^{ij} e_i
\otimes e_j \in {\eu G} \wedge {\eu G} \subset A \wedge A$ gives
now the following Poisson Lie bracket on $G$:
\begin{equation}\label{NilpotBrack}
\{x^p,x^q\} = 2(r^{pl}c_{li}^{q} + r^{lq}c_{li}^{p})x^{i}
\end{equation}
where $\{x^i\} \in {\eu G}^*$ is a basis dual to $\{e_i\}$, and
$c_{ij}^k$ are structure constants of the associative operation
$*$ on $\eu G$. Note that the bracket is nonzero (unlike bracket
on $\eu G$ itself) and linear which suggests possibility to
obtain its quantization.
\end{example}

\end{document}